\documentclass[12pt]{article}
\usepackage{amsmath}
\usepackage{amssymb}

\usepackage[cp1251]{inputenc}
\usepackage[russian]{babel}
\newcommand{\il}[2]{\int\limits_{#1}^{#2}}
\newcommand{\ilp}[1]{\int\limits_{#1}^{+\infty}}

\newcommand{\ph}{\phantom{a}}
\newcommand{\phh}{\phantom{aaa}}

\newcommand{\sist}[2]{\left\{
\begin{array}{l}
{#1}\\
\ph\\
{#2}
\end{array}
\right.}

\begin{document}

MSC 47B35

\vskip 15pt

\centerline{\bf Solvability conditions for a class  of Wiener-Hopf}
\centerline{\bf integral equations of 1-st kind}
\vskip 10pt
\centerline{\bf G. A. Grigorian}

\centerline{\it Institute  of Mathematics of NAS of Armenia}
\centerline{\it E -mail: mathphys2@instmath.sci.am}
\vskip 20 pt

Abstract. The Wiener-Hopf integral equations  of 1-st kind relates to the class of Wiener-Hopf equations of non normal type, to which  the classical Wiener-Hopf method is not applicable, but is completely applicable the special factorization method. In this paper we use the special factorization method to obtain solvability conditions for Wiener-Hopf equations of 1-st kind.

\vskip 10pt

Key words: Wiener-Hopf equation of 1-st kind, the Wiener algebra,  the symbol, the zeroes of the symbol, index of the symbol, index of the operator.

\vskip 10pt

{\bf 1. Introduction}. Let $K(t)$ be an absolutely integrable function on $\mathbb{R}$, and let $f(t)$ be a function on $\mathbb{R}$. Consider the Wiener-Hopf  integral equation  of 1-st kind
$$
 \ilp{0}K(t - \tau) \phi(\tau) d \tau = f(t), \ph t \ge 0.
 \eqno (1.1)
$$
This equation appears in the study of many problems of engineering and physics, in particular,  in the problem of  diffraction of an electromagnetic wave by a perfectly \linebreak conducting half plane (the Sommerfeld problem)   (see [1]). It relates to the class of Wiener-Hopf equations of non normal type, since its symbol degenerates. One of the effective methods for investigating this type of equations is the special factorization method, developed in the book [4]. This method has been used for investigation of Wiener-Hopf equation of second kind in the so called conservative and supercritical cases in the articles [2,3]. In this paper we use the special factorization method to obtain solvability conditions for Eq. (1.1).
Hereafter,
we will assume that $K(t)$ is always representable in the form
$$
K(t) = \sist{\il{t}{\infty}K_1(\tau)d\tau, \ph t > 0,}{\il{-\infty}{t}K_1(\tau) d \tau, \ph t <0} \eqno (1.2)
$$

\vskip 10pt

{\bf 2. Auxiliary propositions}.
 For any $K\in L^1(\mathbb{R})$ we set \ph
$
\nu_m(K) \equiv \il{-\infty}{\infty} t^m K(t) d t, \ph m= 0, 1, 2, \dots.
$
Denote by $E$ one of the spaces $L^p(\mathbb{R}) \ph (1 \le p < \infty), \ph \mathbb{C}^0 \subset \mathbb{M}$, where $\mathbb{M}$ is the space of measurable essentially bounded on $\mathbb{R}$ functions, $\mathbb{C}^0$ is the space of continuous on $\mathbb{R}$ functions $g(t)$ with $\lim\limits_{t \to \pm\infty} g(t) = 0$. We denote by $E_+$ the restriction of $E$ on $[0,+\infty)$, i. e.
$$
E_+ = \{P\phi, \ph \phi \in E\},
$$
where $P$ is the following orthogonal projector:
$$
(P\phi)(t) = \sist{\phi(t), \ph t > 0,}{0, \ph t <0.}
$$
Let $\alpha$ be a real number. We set:
$$
(B_\alpha \phi)(t) \equiv \phi(t) - (1+ i \alpha) e^t\il{t}{\infty}e^{-s}\phi(s) d s, \ph \phi \in E_+, \ph t \ge 0,
$$
$$
(B_\infty)(t) \equiv i e^t \il{t}{\infty} e^{- s} \phi(s) d s, \ph \phi \in E_+, \ph t \ge 0,
$$
$$
(G_\alpha \phi)(t)\equiv \phi(t) + (1 - i\alpha) e^{-i \alpha t} \il{0}{t} e^{i \alpha s} \phi(s) d s, \phi \in E_+, \ph t \ge 0,
$$
$$
(G_\infty f)(t) \equiv i [f(t) + f'(t)], \ph t \ge t_0,
$$
where $f'(t)$ is the generalized derivative of $F(t)$.
Let $\alpha_1, \ldots, \alpha_r$ be any real numbers,  $m_1, \ldots, m_r, m_\infty$ be any natural numbers   (not all necessarily different),    and let $m_1', \ldots, m_r', \linebreak m_\infty, \ph m_1'', \ldots, m_r'', \ph m_\infty''$  be nonnegative integers, such that
$$
m_1' + m_1'' = m_1, \ldots m_r' + m_r'' = m_r, \ph m_\infty' + m_\infty'' = m_\infty.
$$

We put
$$
B\equiv B_\infty^{m_\infty''}\prod_{i=1}^r B_{\alpha_j}^{m_\infty''}, \phh G\equiv G_\infty^{m_\infty''}\prod_{i=1}^r G_{\alpha_j}^{m_\infty''},
$$
$$
\rho_+(\lambda) \equiv (\lambda + i)^{-m_\infty'}\prod\limits_{j=1}^r \biggl(\frac{\lambda - \alpha_j}{\lambda + i}\biggr)^{m_j'}, \ph \rho_-(\lambda) \equiv (\lambda - i)^{-m_\infty''}\prod\limits_{j=1}^r \biggl(\frac{\lambda - \alpha_j}{\lambda - i}\biggr)^{m_j''},
$$
$-\infty \le \lambda \le \infty$. The functions $\rho+(\lambda)$ and $\rho_-(\lambda)$ are called the symbols of the operators $G$ and $B$ respectively (see [4], p. 175) and the spaces
$$
\widetilde{E}_+(\rho_+) \equiv G(E_+), \phh \overline{E}_+(\rho_-) \equiv B(E_+)
$$
are called the spaces, generated by the zeroes of the symbols  $\rho_+$ and $\rho_-$ respectively. The norm $||f||_{\overline{E}_-}$ in $\overline{E}_+(\rho_-)$ is defined as follows
$$
||f||_{\overline{E}_-} = ||B^{-1} f||_{E_+} \ph \mbox{for} \ph  E_+ = L^P_+ (1 \le p < \infty) \ph \mbox{or} \ph E_+ = C_+^0
$$
and
$$
||f||_{\overline{E}_-} = \inf\limits_{B\phi = f} || \phi||_{E_+}, \ph \mbox{for} \ph E_+ = M_+.
$$
By similar way is defined the norm $||f||_{\widetilde{E}_+}$ in $\widetilde{E}_+(\rho_+)$.

\vskip 10pt

Let $k\in L^1(\mathbb{R}),$ and let $c$ be a complex number. Consider the Wiener-Hopf operator
$$
(\widehat{A}\phi) \equiv c \phi(t) - \ilp{0} k(t-s) \phi(s) d s, \phh t \ge 0 \eqno (2.1)
$$
and the corresponding symbol
$$
\mathcal{A}(\lambda) \equiv c - \il{-\infty}{+\infty} e^{ i \lambda t} k(t) d t \ph (-\infty \le \lambda \le +\infty).
$$
Assume
$$
\mathcal{A}(\lambda) = \rho_+(\lambda) \mathcal{C}(\lambda) \rho_-(\lambda), \phh \mathcal{C}(\lambda) \in \mathbb{W},  \phh -\infty \le \lambda \le  +\infty,
$$
where $\mathbb{W}$ is the extended Wiener algebra.

{\bf Theorem 2.1} ([4, p. 175, Theorem 2.6]). {\it  For the operator $\widehat{A}$, defined by (2.1) to be a $\Phi_+$-operator or a $\Phi_-$-operator from $\widetilde{E}_+(\rho_+)$ into $\overline{E}(\rho_-)$ it is necessary and sufficient that $\mathcal{C}(\lambda) \ne 0, \ph -\infty \le \lambda \le +\infty.$ If this condition is satisfied, then the invertibility of the operator $\widehat{A}$ corresponds to the number
$\kappa - \delta$, and the formulae
$$
dimker \widehat{A} = \max \{\delta -\kappa, 0\}, \phh dimcoker \widehat{A} = \max\{\kappa - \delta, 0\}
$$
hold. Here $\kappa = ind \mathcal{C}(\lambda),$
$$
\delta = \sist{0, \ph \mbox{for} \ph E_+ = L_+^p \ph (1 \le p <\infty) \ph \mbox{or} \ph E_+ = C_+^0,}{r'', \ph \mbox{for} \ph E_+ = M_+,}
$$
and $r''$ is the number of positive $m_j'' (j=1,2,\ldots,r).$ In the case $E_+ = C_+$ we have
$$
\delta = \sist{1, \ph \mbox{if} \ph \alpha_{j_0} = 0 \ph \mbox{and} \ph m_{j_0}'' > (1 \le j_0 \le r),}{0, \ph \mbox{otherwise}.}
$$
}

\phantom{aaaaaaaaaaaaaaaaaaaaaaaaaaaaaaaaaaaaaaaaaaaaaaaaaaaaaaaaaaaaaaaaaaaaaa} $\blacksquare$

For any function $k(t)$ on $\mathbb{R}$ we set $\nu_j(k)\equiv \il{-\infty}{\infty} t^j k(t) d t, \ph j=0,1,2, \ldots .$
Consider the Wiener-Hopf operator
$$
(\widehat{a}\phi)(t) \equiv \ilp{0} K(t -\tau) \phi(\tau) d \tau
$$
and the corresponding symbol
$$
a(\lambda) \equiv \il{-\infty}{\infty} e^{i \lambda t} K(t) d t, \phh -\infty \le \lambda \le \infty.
$$
Assume $\nu_1(|K_1|) < \infty$. Then from (1.2) we obtain
$$
a(\lambda) = \il{0}{\infty} e^{i\lambda t} d t \il{t}{\infty} K_1(\tau) d \tau + \il{-\infty}{0} e^{i\lambda t} d t\il{-\infty}{t} K_1(\tau) d \tau =
$$
$$
= \il{0}{\infty}K_1(\tau) d \tau \il{0}{\tau} e^{i\lambda t} d t + \il{-\infty}{0} K_1(\tau) d \tau \il{\tau}{0} e^{i\lambda t} d t =
$$
$$
=\frac{i}{\lambda} \biggl\{\nu_0(\widetilde{K}_1) - \il{-\infty}{\infty} e^{i\lambda t} \widetilde{K}_1(t) d t\biggr\}, \ph -\infty \le \lambda \le \infty, \ph \lambda \ne 0, \eqno (2.2)
$$
where $\widetilde{K}_1(t) \equiv \sist{K_1(t), \ph t > 0,}{-K_1(t), \ph t < 0.}$ It is not difficult to verify that
$$
b(\lambda) \equiv \nu_0(\widetilde{K}_1) - \il{-\infty}{\infty} e^{i\lambda t} \widetilde{K}_1(t) d t =
$$
$$
= \il{0}{\infty}(1- \cos \lambda t)[K_1(t) - K_1(-t)] dt - i \il{-\infty}{\infty} \sin \lambda t \widetilde{K}_1(t) d t, \ph -\infty \le \lambda \le \infty. \eqno (2.3)
$$
Assume
$$
K_1(t) - K_1(-t) \ge 0, \ph t \ge t \ge 0 \eqno (2.4)
$$
and
$$
\il{0}{\infty}[K_1(t) - K_1(-t)] d t > 0. \eqno (2.5)
$$
Then it follows from (2.2) and (2.3)  that
$$
a(\lambda) = \frac{i}{\lambda} b(\lambda), \ph b(\lambda) \ne 0, \ph -\infty \le \lambda \le \infty, \ph \lambda \ne 0
$$
and
$$
b(0) = \nu_0(\widetilde{K}_1) - \nu_0(\widetilde{K}_1) = 0. \eqno (2.6)
$$
It follows from (2.4) and (2.5) that
$$
b(\pm\infty) = \nu_0(\widetilde{K}_1) > 0. \eqno (2.7)
$$
The relations (2.3)-(2.5) imply that
$$
-\frac{\pi}{2} < arg \hskip 3pt b(\lambda) < \frac{\pi}{2}, \ph -\infty < \lambda < \infty, \ph \lambda \ne 0. \eqno (2.8)
$$
Let us discuss the properties of $a(\lambda)$ for the cases  $\nu_1(\widetilde{K}_1) >0, \ph \nu_1(\widetilde{K}_1) <0, \ph\nu_1(\widetilde{K}_1) = 0.$

$I. \ph  \nu_1(\widetilde{K}_1) >0.$ Consider the function
$$
c(\lambda) \equiv \Bigl[1 + \frac{i}{\lambda}\Bigr] b(\lambda), \ph \lambda \ne 0.
$$
It follows from (2.5) and (2.7) that
$$
c(\pm\infty) = b(\pm\infty) = \nu_0(\widetilde{K}_0) > 0. \eqno(2.9)
$$
We have
$$
c(0) = \lim\limits_{\lambda \to 0} b(\lambda) + i\lim\limits_{\lambda \to 0}\frac{b(\lambda)}{\lambda} =  i \lim\limits_{\lambda \to 0}\ilp{0}\frac{(1 - \cos \lambda t)}{\lambda} [K_1(t) - K_1(-t)] d t +
$$
$$
+\lim\limits_{\lambda \to 0}\il{-\infty}{\infty}\frac{\sin \lambda t}{\lambda} \widetilde{K}_1(t) d t = \nu_1(\widetilde{K}_1) > 0. \eqno (2.10)
$$
Obviously $-\frac{\pi}{2} < 1 + \frac{i}{\lambda} < \frac{\pi}{2}, \ph \lambda \ne 0$. This together with (2.8)--(2.10) implies that \linebreak $-\pi < arg \hskip 3pt c(\lambda) < \pi, \ph -\infty \le \lambda \le +\infty.$ Hence,
$$
a(\lambda) =  \frac{1}{\lambda + i} c_1(\lambda), \ph c_1(\lambda) \in \mathbb{W}, \ph c_1(\lambda) \ne 0, \ph -\infty \le \lambda \le +\infty. \eqno (2.11)
$$
$$
ind \hskip 3pt c_1(\lambda) = 0, \eqno (2.12)
$$
where $c_1(\lambda) = i c(\lambda).$

II. $\nu_1(\widetilde{K}) < 0.$ For this case we have
$$
c(\lambda) = \Bigl(\frac{\lambda + i}{\lambda}\Bigr) b(\lambda) = \Bigl(\frac{\lambda + i}{\lambda - i}\Bigr) \Bigl(1 - \frac{i}{\lambda}\Bigr) b(\lambda).
$$
Hence,
$$
\widetilde{c}(\lambda) \equiv \Bigl(\frac{\lambda - i}{\lambda + i}\Bigr) c(\lambda) = \Bigr[1 - \frac{i}{\lambda}\Bigr] b(\lambda),
$$
$$
\widetilde{c}(\pm\infty) =  b(\pm\infty) = \nu_0(\widetilde{K}_1) > 0, \eqno (2.13)
$$
$$
\widetilde{c}(0) = - i \lim\limits_{\lambda \to 0}\ilp{0}\frac{(1 - \cos \lambda t)}{\lambda} [k_1(t) - k_1(-t)] d t - \lim\limits_{\lambda \to 0} \il{-\infty}{\infty} \frac{\sin \lambda t}{\lambda} \widetilde{K}_1(t) d t = -\nu_1(\widetilde{K}_1) > 0. \eqno (2.14)
$$
Obviously, $-\frac{\pi}{2} < 1 - \frac{i}{\lambda} < \frac{\pi}{2}, \ph \lambda \ne 0$. This together with (2.8), (2.13) and (2.14) implies that
$$
-\frac{\pi}{2} < arg \hskip 3pt \widetilde{c}(\lambda) < \frac{\pi}{2}, \ph -\infty \le \lambda \le +\infty.
$$
Hence, $\widetilde{c}(\lambda) \ne 0, \ph -\infty \le \lambda \le \infty, \ph  ind \hskip 3pt \widetilde{c}(\lambda) = 0.$ From here it follows
$$
a(\lambda) = \frac{1}{\lambda + i} c_1(\lambda), \ph c_1(\lambda) \in \mathbb{W}, \ph c_1(\lambda) \ne 0, \ph -\infty \le \lambda \le \infty, \ph ind \hskip 3pt c_1(\lambda) = -1, \eqno (2.15)
$$
where $c_1(\lambda) = i c(\lambda), \ph -\infty \le \lambda \le \infty.$

III. \ph $\nu_1(\widetilde{K}_1) = 0, \ph \nu_2(\widetilde{K}_1) < \infty \ph (\nu_2(\widetilde{K}_1) = \ilp{0}t^2[K_1(t) - K_1(-t)] d t > 0)$. We have
$$
\lim\limits_{\lambda \to 0}\frac{b(\lambda)}{\lambda^2} = \lim\limits_{\lambda \to 0} \ilp{0}\frac{1 - \cos \lambda t}{\lambda^2}[K_1(t) - K_1(-t)] d t - i \il{-\infty}{\infty} \frac{\sin \lambda t}{\lambda^2} K_1(t) d t = \frac{1}{2} \nu_2(\widetilde{K}_1) > 0. \eqno (2.16)
$$
Therefore ([4, p. 180, Theorem 2.9]),
$$
b(\lambda) = \frac{\lambda^2}{1 + \lambda^2} c(\lambda), \ph c(\lambda) \in \mathbb{W}, \ph -\infty \le \lambda \le \infty, \eqno (2.17)
$$
We have $c(\lambda) = [1 + \frac{1}{\lambda^2}] b(\lambda)$ and
$$
c(\pm\infty) = b(\pm\infty) = \nu_0(\widetilde{K}_1) > 0. \eqno (2.18)
$$
By (2.16) we have $c(0) = \frac{1}{2} \nu_2(\widetilde{K}_1) > 0$. This together with (2.8) and (2.18) implies that
$$
-\frac{\pi}{2} < arg \hskip 3pt c(\lambda) < \frac{\pi}{2}, \phh -\infty \le \lambda \le \infty.
$$
Hence,
$$
a(\lambda) = \frac{\lambda}{1 + \lambda^2} c_1(\lambda), \ph c_1(\lambda) \in \mathbb{W}, \ph c_1(\lambda) \ne 0, \ph -\infty \le \lambda \le \infty, \ph ind \hskip 3pt c_1(\lambda) = 0, \eqno (2.19)
$$
where $c_1(\lambda) = i c(\lambda), \ph -\infty \le \lambda \le \infty.$

Now we consider the case

\vskip 10pt

\noindent
$
(\alpha) \ph K_1(t) = \sist{\ilp{t} K_0(\tau) d \tau, \ph t > 0,}{\il{-\infty}{t} K_0(\tau) d \tau, \ph t < 0,} \ph K_0(t) + K_0(-t) \ge 0, \ph t \ge 0, \ph \nu_j(|K_0|) < \infty, \linebreak j=0,1,2,$

\vskip 10pt

\noindent
$
(\beta)
\ph \nu_1(K_0) = 0, \ph \nu_0(K_0) > 0, \ph \nu_2(K_0) > 0.
$

\vskip 10pt

\noindent
Then by (2.2)
$$
a(\lambda) = \frac{i}{\lambda} \biggl\{\ilp{0}d\tau\il{\tau}{+\infty} K_0(u) d u - \il{-\infty}{0}d\tau\il{-\infty}{\tau} K_0(u) d u - \ilp{0}e^{i\lambda \tau} d \tau \il{\tau}{+\infty} K_0(u) d u +
$$
$$
+ \il{-\infty}{0}e^{i \lambda \tau} d\tau\il{-\infty}{\tau} k_0(u) d u\biggr\} = \frac{i}{\lambda}\biggl\{\nu_1(K_0) - \frac{1}{i\lambda}\biggl[\ilp{0}(e^{i\lambda u} - 1) K_0(u) d u + \il{-\infty}{0}(e^{i\lambda u} - 1) K_0(u) d u \biggr] \biggr\}=
$$
$$
= \frac{1}{\lambda^2}\biggl[\nu_0(K_0) - \il{-\infty}{\infty}e^{i \lambda u} K_0(u) d u\biggr], \phh -\infty \le \lambda \le +\infty.  \eqno (2.20)
$$
Thus
$$
a(\lambda) = \frac{1}{\lambda^2} d(\lambda), \phh d(\lambda) \in \mathbb{W}, \eqno (2.21)
$$
where
$d(\lambda) \equiv \nu_0(K_0) - \il{-\infty}{\infty}e^{i \lambda u} K_0(u) d u, \ph -\infty \le \lambda \le +\infty.$
We have
$$
d(\pm\infty) = \nu_0(K_0) > 0, \ph d(0) = \nu_0(K_0) -  \nu_0(K_0)  = 0, \eqno (2.22)
$$
$$
d \hskip 2pt '(0) = - i\il{-\infty}{\infty}t K_0(t) d t = 0, \ph d \hskip 2pt ''(0) = i\il{-\infty}{\infty}t^2 K_0(t) d t \ne 0.
$$
Therefore
$$
d(\lambda) = \frac{\lambda^2}{1 + \lambda^2} e(\lambda) \phh e(\lambda) \in \mathbb{W}, \eqno (2.23)
$$
$$
e(\pm\infty) = d(\pm\infty) = \nu_0(K_0) > 0, \ph
e(0) = \lim\limits_{\lambda \to 0}\Bigl[1 + \frac{1}{\lambda^2}\Bigr] d(\lambda) =
$$
$$
=\lim\limits_{\lambda \to 0} \il{-\infty}{\infty}\frac{(1 - \cos \lambda u)}{\lambda^2} K_0(u) d u - i \il{-\infty}{\infty}\frac{\sin \lambda u}{\lambda^2} K_0(u) d u = \frac{1}{2}\nu_2(K_0) -
$$
$$
- i \lim\limits_{\lambda \to 0} \il{-\infty}{\infty}\frac{\cos \lambda u}{2\lambda} u K_0(u) d u = \frac{1}{2}\nu_2(K_0) + i \lim\limits_{\lambda \to 0} \il{-\infty}{\infty}\frac{(1 - \cos \lambda u)}{2\lambda} u K_0(u) d u =
$$
$$
= \frac{1}{2}\nu_2(K_0) +
 i \lim\limits_{\lambda \to 0} \il{-\infty}{\infty}\frac{\sin \lambda u}{2\lambda} u^2 K_0(u) d u = \frac{1}{2}\nu_2(K_0) > 0, \eqno (2.24)
$$
This together with (2.21), (2.23) implies
$$
a(\lambda) = \frac{1}{1 +\lambda^2} e(\lambda), \phh e(\lambda) \in \mathbb{W}, \phh e(\lambda) \ne 0, \phh -\infty\le \lambda \le +\infty. \eqno (2.25)
$$
It follows from here  from the conditions $a_0(t) + a_0(-t) \ge 0, \ph t\ge t_0, \ph \nu_0(K_0) > 0$ and   from (2.20),  (2.24) that
$$
ind \hskip 3pt e(\lambda) = 0. \eqno (2.26)
$$

\vskip 10pt

{\bf 3. Main result.} According to Theorem 2.1 it follows from  (2.11) and (2.12) that  if $\nu_1(\widetilde{K}_1) > 0$, then  for every $f(t) \in E_+$ Eq. (1.1) has a solution in $\widetilde{E}_+(\frac{1}{\lambda + i})$. The corresponding homogeneous equation  has only the trivial solution on $\widetilde{E}_+(\frac{1}{\lambda + i})$.
 It follows from (2.15) that
$$
a(\lambda) = \frac{1}{\lambda - i}\biggl[\biggl(\frac{\lambda - i}{\lambda + i}\biggr) c_1(\lambda)\biggr], \phh -\infty \le \lambda \le +\infty \phh ind \hskip 3pt \biggl[\biggl(\frac{\lambda - i}{\lambda + i}\biggr) c_1(\lambda)\biggr] = 0.
$$
Then   by Theorem 2.1 if $\nu_1(\widetilde{K}_1)~ <~ 0$, then for every $f(t) \in E_+$ Eq. (1.1) has a solution in   $\widetilde{E}_+(\frac{1}{\lambda + i})$. The corresponding homogeneous equation  has only one linearly independent solution in $\widetilde{E}_+(\frac{1}{\lambda + i})$.
Hence by Theorem 2.1  it follows from (2.15) that if $\nu_1(\widetilde{K}_1) < 0$ and $f(t) \in \overline{E}_+(\frac{1}{\lambda - i})$, then the solution $\phi(t)$ of Eq. (1.1) belongs to $E_+$.
If $\nu_1(\widetilde{K}_1) =~ 0, \ph \nu_2(\widetilde{K}_1) <~ \infty \ph (\nu_2(\widetilde{K}_1) = \ilp{0}t^2[K_1(t) - K_1(-t)] d t > 0)$, then it follows from (2.19) that
$$
a(\lambda) = \frac{\lambda}{(\lambda + i)^2}\biggl[\biggl(\frac{\lambda + i}{\lambda - i}\biggr) c_1(\lambda)\biggr], \phh -\infty \le \lambda \le +\infty, \phh ind \hskip 3pt \biggl[\biggl(\frac{\lambda - i}{\lambda + i}\biggr) c_1(\lambda)\biggr] = -1.
$$
By Theorem 2.1 it follows from here and from (2.19) that for every $f(t)\in E_+$ Eq. (1.1) has a solution $\phi(t) \in \widetilde{E}(\frac{\lambda}{(\lambda + i)^2}),$ the corresponding homogeneous equation has only one linearly independent solution $\phi_0(t)$ in $\widetilde{E}(\frac{\lambda}{(\lambda + i)^2})$. It is clear from (2.19) that $\phi_0(t)$ does not belong to $\widetilde{E}_+(\frac{\lambda}{\lambda + i})$. Therefore,
$\phi_0(t) \in \widetilde{E}(\frac{\lambda}{(\lambda + i)^2}) \backslash \widetilde{E}_+(\frac{\lambda}{\lambda + i})$. If  $f(t) \in \overline{E}_+(\frac{1}{\lambda - i})$, then  $\phi(t) \in \widetilde{E}_+(\frac{\lambda}{\lambda + i})$.
 Let the conditions $(\alpha)$ and $(\beta)$ be satisfied. Then it follows from (2.25) and (2.26) that
 $$
a(\lambda) = \frac{1}{(\lambda + i)^2}\biggl[\biggl(\frac{\lambda + i}{\lambda - i}\biggr) e(\lambda)\biggr], \phh -\infty \le \lambda \le +\infty, \phh ind \hskip 3pt \biggl[\biggl(\frac{\lambda - i}{\lambda + i}\biggr) e (\lambda)\biggr] = -1.
$$
Then by Theorem 2.1 it follows from  here  that  for every $f(t) \in E_+$ Eq. (1.1) has a solution $\phi(t) \in \widetilde{E}(\frac{1}{(\lambda + i)^2}),$ the corresponding homogeneous equation has only one linearly independent solution $\phi_0(t)$ in $\widetilde{E}(\frac{1}{(\lambda + i)^2})$. Moreover according to (2.25) and (2.26) if  $f(t) \in \overline{E}_+(\frac{1}{\lambda - i})$, then  $\phi(t) \in \widetilde{E}_+(\frac{1}{\lambda + i})$ and the corresponding homogeneous equation has only the trivial solution in $\widetilde{E}_+(\frac{1}{\lambda + i})$. Hence,  $\phi_0(t) \in \widetilde{E}(\frac{1}{(\lambda + i)^2}) \backslash \widetilde{E}_+(\frac{1}{\lambda + i})$.

Let us now summarize the results obtained in the form

\vskip 20pt

{\bf Theorem 3.1. \ph (Main result).} {\it Let the kernel function $K(t)$ has the representation of the form (1.2). Then the following assertions are valid.

\noindent
1) \ph If  $\nu_1(\widetilde{K}_1) > 0$, then Eq. (1.1) has a solution in $\widetilde{E}_+(\frac{1}{\lambda + i})$. The corresponding \linebreak homogeneous equation  has  only the trivial solution in $\widetilde{E}_+(\frac{1}{\lambda + i})$.

\noindent
2) \ph If  $\nu_1(\widetilde{K}_1) < 0$, then  Eq. (1.1) has a solution  in $\widetilde{E}_+(\frac{1}{\lambda + i})$. The corresponding homogeneous equation  has only one linearly independent solution in $\widetilde{E}(\frac{\lambda}{(\lambda + i)^2}) \backslash \widetilde{E}_+(\frac{\lambda}{\lambda + i})$. If  $f(t) \in \overline{E}_+(\frac{1}{\lambda - i})$, then  $\phi(t) \in \widetilde{E}_+(\frac{\lambda}{\lambda + i})$.

\noindent
3) \ph If $\nu_1(\widetilde{K}_1) = 0, \ph \nu_2(\widetilde{K}_1) < \infty \ph (\nu_2(\widetilde{K}_1) = \ilp{0}t^2[K_1(t) - K_1(-t)] d t > 0)$, then for every $f(t) \in E_+$ Eq. (1.1) has a solution in  $\widetilde{E}_+(\frac{\lambda}{(\lambda + i)^2})$. The corresponding homogeneous equation has only one linearly independent  solution in $\widetilde{E}_+(\frac{\lambda}{(\lambda + i)^2})$.

\noindent
4) \ph If  the conditions $(\alpha)$ and $(\beta)$ hold, then for every  $f(t) \in E_+$ Eq. (1.1) has a solution $\phi(t) \in \widetilde{E}(\frac{1}{(\lambda + i)^2}),$ the corresponding homogeneous equation has only one linearly \linebreak independent solution in $\widetilde{E}(\frac{1}{(\lambda + i)^2}) \backslash \widetilde{E}_+(\frac{1}{\lambda + i})$. If  $f(t) \in \overline{E}_+(\frac{1}{\lambda - i})$, then  $\phi(t) \in \widetilde{E}_+(\frac{1}{\lambda + i})$.
}

\vskip 10pt

{\bf Remark 3.1.} {\it The main result of work [1] (see [1, Theorem 3.1]) is conditional in the sense that one of its conditions requires that the symbol function $G(\omega)$ (associated with the kernel function $k(t)$)  must not vanish on $\mathbb{R}$ (the condition (II)). Verification  of this condition is  difficult, and makes it hard to use the mentioned  main result to concrete equations. While the conditions of the  result obtained  of this paper are enough verifiable.
}

\vskip 20pt

\centerline{\bf References}

\vskip 20pt

\noindent
1. A. F. Dos Santos and F. S. Teixsera, Theory of a Class of Wiener-Hopf Equations of the  \linebreak \phantom{a} First Kind: Application to the Sommerfeld Problem. J. Math. Anal. Appl.  128,  \linebreak \phantom{a}189-204 (1987)

\noindent
2. G. A. Grigorian, Solvability of a class of Wiener-Hopf integral equations. Izv. Nats. \linebreak \phantom{a} Akad. Nauk. Armenii. Matematika, vol 21, No. 2, 1996, pp. 21-32.

\noindent
3. G. A. Grigorian The Wiener-Hopf equation in the supercritical case.
Izv. Nats.  Akad.  \linebreak \phantom{a} Nauk Armenii Mat. 32 (1997), no. 1, 60–74.

\noindent
4. S. Presdorf, Some Classes of Singular Equations, Mir, 1979.

\end{document}